\documentclass [12pt]{article}
\usepackage{graphicx,amssymb,amsfonts,latexsym,amsmath,amsthm,times}
\usepackage{epsfig}
\usepackage{color}
\usepackage[english]{babel}
\usepackage[a4paper]{geometry}
\def\lanbox{\hbox{$\, \vrule height 0.25cm width 0.25cm depth 0.01cm \,$}}

\numberwithin{equation}{section}

\begin{document}

\vspace*{1.4cm}

\normalsize \centerline{\Large \bf ON THE EXISTENCE OF SOLUTIONS
FOR A CLASS OF}

\medskip

\centerline{\Large\bf SYSTEMS OF INTEGRO-DIFFERENTIAL EQUATIONS}

\medskip

\centerline{\Large\bf WITH THE LOGARITHMIC LAPLACIAN AND DRIFT}

\vspace*{1cm}

\centerline{\bf Yuming Chen$^{1}$, Vitali Vougalter$^{2 \ *}$}

\vspace*{0.5cm}

\centerline{$^1$ Department of Mathematics, Wilfrid Laurier University}

\centerline{Waterloo, Ontario, N2L 3C5, Canada}

\centerline{ e-mail: ychen@wlu.ca}

\medskip

\centerline{$^{2 \ *}$ Department of Mathematics, University
of Toronto}

\centerline{Toronto, Ontario, M5S 2E4, Canada}

\centerline{ e-mail: vitali@math.toronto.edu}

\medskip


\vspace*{0.25cm}

\noindent {\bf Abstract:}
In this article, we consider a system of integro-differential equations in
$L^2(\mathbb{R}, \mathbb{R}^{N})$, which  contains the logarithmic Laplacian in
the presence of transport terms. The linear operators associated with the
system satisfy the Fredholm property. By virtue of a fixed point technique,
we demonstrate the existence of solutions. We emphasize that the discussion
is more complicated than that of the scalar situation as there are
more cumbersome technicalities to overcome.

\vspace*{0.25cm}



\vspace*{0.5cm}

\bigskip

\bigskip


\setcounter{section}{1}

\centerline{\bf 1. Introduction}

\medskip

\noindent
The present work is devoted to the solvability of the following system of
integro-differential equations with $1\leq m\leq N$, 
$$
\Big[-\frac{1}{2}\hbox{ln}\Big(-\frac{d^{2}}{dx^{2}}\Big)\Big]u_{m}+
$$
\begin{equation}
\label{h}
b_{m}\frac{du_{m}}{dx}+a_{m}u_{m}+f_{m}(x)+
\int_{-\infty}^{\infty}K_{m}(x-y)g_{m}(u(y))dy=0, 
\end{equation}
where $x\in \mathbb{R}$ and $a_{m}, b_{m}\in\mathbb{R}, \ b_{m}\neq 0$ are
constants.  Let us recall that the solvability of a single equation resembling
~\eqref{h}  with the fractional Laplacian in the diffusion term has been
considered by Vougalter and Volpert~\cite{VV21}.

Biological processes  
with nonlocal consumption of resources and  intra-specific competition can be
described by nonlocal reaction-diffusion equations (see, for example, \cite{ABVV10,BNPR09,GVA06,VV130} and
references therein).
In such time-dependent problems the space variable $x$ corresponds to the cell
genotype and $u_{m}(x,t)$ stand for the cell density distributions for various
groups of cells as functions of their genotype and time,
\begin{equation}
\label{vf}
u(x,t)=(u_{1}(x,t), u_{2}(x,t), ..., u_{N}(x,t))^{T}.  
\end{equation}  
The evolution of the cell densities occurs via 
cell proliferation, mutations, transport, and cell influx/efflux.
The diffusion terms correspond to the change of genotype
due to small random mutations  and the integral terms describe large mutations.
The functions $g_{m}(u)$ denote the rates of cell births depending on $u$
(density-dependent proliferation). The kernels $K_{m}(x-y)$ give
the proportions of newly born cells changing their genotypes from $y$ to $x$.
We assume that they depend on the distances between the genotypes.
The source terms $f_{m}(x)$ designate
the influxes/effluxes of cells for different genotypes.
Since models of this kind describe the distributions of population
densities with respect to the genotype, the existence of stationary solutions
of the nonlocal reaction-diffusion problems
corresponds to the existence of biological species.

We introduce the concept of  non-Fredholm operators. Consider the
problem
\begin{equation}
\label{eq1}
 -\Delta u + V(x) u - a u=f,
\end{equation}
where $u \in E= H^{2}({\mathbb R}^{d})$ and  $f \in F=
L^{2}({\mathbb R}^{d})$, $d\in {\mathbb N}$, $a$ is a constant, and
the scalar potential function $V(x)$ either vanishes identically
or converges to $0$ at infinity. If $a \geq 0$, the origin belongs to the
essential spectrum of the operator $A : E \to F$ corresponding to the left
side of equation ~\eqref{eq1}. Consequently, such an operator fails to
satisfy the Fredholm property. Its image is not closed, and, for $d>1$,
the dimension of its kernel and the codimension of its image are
not finite.
Note that for the Fredholm property we consider bounded but not
$H^{2}({\mathbb R}^{d})$ solutions of the corresponding homogeneous adjoint
equation.
If $V(x)$ vanishes identically, problem~\eqref{eq1} has
constant coefficients and we can apply the Fourier transform to solve it
explicitly. If $f\in L^{2}({\mathbb R}^{d})$ and $xf\in L^{1}({\mathbb R}^{d})$,
then it has a unique solution in $H^{2}({\mathbb R}^{d})$ if and only if
\begin{equation}
\label{ocb}
\Bigg(f(x), \frac{e^{ipx}}{(2\pi)^{\frac{d}{2}}}\Bigg)_{L^{2}({\mathbb R}^{d})}=0, \quad
p\in S_{\sqrt{a}}^{d} \quad a.e.
\end{equation}
(see Lemmas 5 and 6 of ~\cite{VV14}). Here $S_{\sqrt{a}}^{d}$ stands for the
sphere in ${\mathbb R}^{d}$ of radius $\sqrt{a}$ centered at the origin.
Hence,  though our operator does not satisfy the Fredholm property, the
solvability conditions can be formulated similarly. But this similarity is only formal since the range of the operator is not closed. The orthogonality
conditions~\eqref{ocb}  are with respect to the standard Fourier harmonics,
which solve the homogeneous adjoint problem for~\eqref{eq1}  when the
scalar potential function is trivial. They
belong to $L^{\infty}({\mathbb R}^{d})$ but they are not square integrable.

Let us mention  that elliptic problems containing
non-Fredholm operators have been studied actively in recent years.
Approaches in weighted Sobolev and H\"older spaces have been developed in~\cite{Amrouche1997,Amrouche2008,Bolley1993,Bolley2001,B88}. Particularly,
when $a=0$, the operator $A$ is Fredholm in certain properly chosen
weighted
spaces (see~\cite{Amrouche1997,Amrouche2008,Bolley1993,Bolley2001,B88}).
However, the situation when $a \neq 0$ is considerably
different and the approaches developed in these articles are not applicable.
Non-Fredholm Schr\"odinger type operators
were treated with the methods of the spectral and the
scattering theory in~\cite{EV21,V2011,VV08,VV19}.
Fredholm structures, topological invariants, and applications were discussed
in~\cite{E09}. The articles ~\cite{GS05} and~\cite{RS01} deal with the
Fredholm and properness properties of  quasilinear elliptic
systems of the second order and of the operators of this kind on
${\mathbb R}^{N}$. The exponential decay and Fredholm properties in  
second-order quasilinear elliptic systems of equations were covered in
~\cite{GS10}.
The nonlinear non-Fredholm elliptic equations were discussed in
~\cite{EV20,VV14,VV21}. Considerable applications to the theory of
reaction-diffusion
equations were developed in ~\cite{DMV05,DMV08}.

The logarithmic Laplacian $\hbox{ln}(-\Delta)$ is the operator with Fourier
symbol $2\ln |p|$. It appears as the formal derivative
$\partial_{s}|_{s=0}(-\Delta)^{s}$ of fractional Laplacians at $s=0$.
The operator $(-\Delta)^{s}$ has been extensively  used, for instance, in the studies of  anomalous diffusion problems (to name a few, see~\cite{MK00,VV19,VV21} and  references therein). Spectral properties of the
logarithmic Laplacian in an open set of finite measure with Dirichlet boundary
conditions were investigated in~\cite{LW21} (see also ~\cite{CW19}).
The studies of $\ln (-\Delta)$ are crucial to understand the asymptotic spectral properties of the family of fractional Laplacians
in the limit $s\to 0^{+}$. In~\cite{JSW20}, it was established that this
operator enables us to characterize the $s$-dependence of  solutions to
fractional Poisson problems for the full range of exponents $s\in (0, 1)$.

The operators involved in the left side of system ~\eqref{h} are given by
\begin{equation}
\label{lab}
L_{a_{m}, b_{m}}:=\frac{1}{2}\ln\Big(-\frac{d^{2}}{dx^{2}}\Big)-b_{m}
\frac{d}{dx}-a_{m},
\quad a_{m},b_{m}\in {\mathbb R}, \quad b_{m}\neq 0
\end{equation}
with $1\leq m\leq N$
and are considered on $L^{2}({\mathbb R})$. By means of the standard Fourier
transform, it can be easily derived that the essential spectrum of
(\ref{lab}) is equal to
\begin{equation}
\label{ess}
\lambda_{a_{m}, b_{m}}(p)=\ln \Big(\frac{|p|}{e^{a_{m}}}\Big)-ib_{m}p, \quad
a_{m},b_{m}\in {\mathbb R}, \quad b_{m}\neq 0, \quad 1\leq m\leq N.
\end{equation}
Evidently, for $p\in {\mathbb R}, \ 1\leq m\leq N$,
\begin{equation}
\label{lbab}
|\lambda_{a_{m},  b_{m}}(p)|=\sqrt{\ln^{2}\Big(\frac{|p|}{e^{a_{m}}}\Big)+
b_{m}^{2}p^{2}}\geq C>0,
\end{equation}
where $C$ is a constant. Let us note that as distinct from the case without
the drift term studied in~\cite{EV23}, each operator~\eqref{lab}  has
the Fredholm property. Solvability of certain nonhomogeneous linear problems
involving the logarithmic Schr\"odinger operators in higher dimensions was
considered in~\cite{EV231}. The logarithmic Schr\"odinger operator and
associated Dirichlet problems were discussed in ~\cite{F23}. The logarithmic
Laplacian also arises in the geometric context of the $0$-fractional perimeter,
which has been studied in ~\cite{DNP21}. The article ~\cite{ZKR23} deals with
the symmetry of positive solutions for Lane-Emden systems involving the
logarithmic Laplacian.

We set $K_{m}(x) = \varepsilon_{m} {\cal K}_{m}(x), \ 1\leq m\leq N$,
where all $\varepsilon_{m} \geq 0$ and introduce
\begin{equation}
\label{epsm}
\varepsilon:=\hbox{max}_{1\leq m\leq N}\varepsilon_{m}.
\end{equation}
In the biological applications, such small and nonnegative parameters
$\varepsilon_{m}$ have the meaning
that the integral production terms are small with respect to the others, that
is, the frequency of large mutations is sufficiently small.
Our first assumption is as follows.

\noindent {\bf Assumption 1.}  {\it Let $1\leq m\leq N$, the constants
$a_{m}$, $b_{m}\in {\mathbb R}$ with $b_{m}\neq 0$.
The functions $f_{m}(x): {\mathbb R}\to {\mathbb R}$ do not vanish identically
on the real line for some $1\leq m\leq N$ and $f_{m}(x)\in L^{2}({\mathbb R})$.
Let us also assume that ${\cal K}_{m}(x): {\mathbb R}\to {\mathbb R}$ and 
${\cal K}_{m}(x)\in L^{1}({\mathbb R})$. Moreover,}
\begin{equation}
\label{K}
{\cal K}^{2}:=\sum_{m=1}^{N}\|{\cal K}_{m}\|_{L^{1}({\mathbb R})}^{2}>0.  
\end{equation}

For a vector function
\begin{equation}
\label{vfst}
u(x)=(u_{1}(x), u_{2}(x), ..., u_{N}(x))^{T},
\end{equation}  
we will use the norm
\begin{equation}
\label{nvf}
\|u\|_{L^{2}({\mathbb R}, {\mathbb R}^{N})}^{2}:=\sum_{m=1}^{N}\|u_{m}\|_{L^{2}({\mathbb R})}^{2}.
\end{equation}
  
When all the nonnegative parameters $\varepsilon_{m}$ vanish, we obtain the
linear system of equations~\eqref{lp}.
By virtue of Lemma 5 further down along with Assumption 1 and Remark 6,
system ~\eqref{lp} has a unique solution
\begin{equation}
\label{u0}
u_{0}(x)\in L^{2}({\mathbb R}, {\mathbb R}^{N}),
\end{equation}
which does not vanish identically on the real line.
Let us seek the resulting solution of the nonlinear system of equations
~\eqref{h}  as
\begin{equation}
\label{r}
u(x)=u_{0}(x)+u_{p}(x),
\end{equation}
where
\begin{equation}
\label{up}
u_{p}(x):=(u_{p, 1}(x), u_{p, 2}(x), ..., u_{p, N}(x))^{T}. 
\end{equation}  
Obviously, we arrive at the perturbed system
$$
\Big[\frac{1}{2}\ln \Big(-\frac{d^{2}}{d x^{2}}\Big)\Big]u_{p, m}-
b_{m}\frac{du_{p, m}}{dx}-a_{m}u_{p, m}
=
$$
\begin{equation}
\label{pert}
\varepsilon_{m} \int_{-\infty}^{\infty}
{\cal K}_{m}(x-y)g_{m}(u_{0}(y)+u_{p}(y))dy
\end{equation}
with $1\leq m\leq N$.
For technical purpose, we introduce a closed ball in our function space,
namely,
\begin{equation}
\label{b}
B_{\rho}:=\{u(x)\in L^{2}({\mathbb R}, {\mathbb R}^{N}) \ | \
\|u\|_{L^{2}({\mathbb R}, {\mathbb R}^{N})}\leq
\rho \}, \quad 0<\rho\leq 1.
\end{equation}
Let us look for solutions of the system of equations ~\eqref{pert} as  fixed
points of the auxiliary nonlinear problem,
$$
\Big[\frac{1}{2}\ln \Big(-\frac{d^{2}}{d x^{2}}\Big)\Big]u_{m}-
b_{m}\frac{du_{m}}{dx}-a_{m}u_{m}
=
$$
\begin{equation}
\label{aux}
\varepsilon_{m} \int_{-\infty}^{\infty}
{\cal K}_{m}(x-y)g_{m}(u_{0}(y)+v(y))dy,
\end{equation}
where $1\leq m\leq N$,
in the ball $B_{\rho}$ given by ~\eqref{b}. The fixed point technique was
used in ~\cite{CV21} to demonstrate the
persistence of pulses for a class of reaction-diffusion equations. For a given
vector function $v(y)$, system \eqref{aux} is to be solved for $u(x)$. 

The left side of ~\eqref{aux} involves the Fredholm operators given by
~\eqref{lab}. The solvability of the linear and nonhomogeneous system
containing these operators will be considered in the final section of this
article. Let us point out that, as distinct from the present case, the
equations discussed in ~\cite{EV20,VV14} involved
operators without the Fredholm property and orthogonality
relations were used to solve them.

For our nonlinear vector function
\begin{equation}
\label{g}
g(z):=(g_{1}(z), g_{2}(z), ..., g_{N}(z)),  
\end{equation}
its gradient is given by
\begin{equation}
\label{grg}
\nabla g(z):=(\nabla g_{1}(z), \nabla g_{2}(z), ..., \nabla g_{N}(z)),  
\end{equation}
such that the norm
\begin{equation}
\label{ng}
\|\nabla g(z)\|_{L^{\infty}({\mathbb R}^{N}, {\mathbb R}^{N^{2}})}:=
\sum_{m=1}^{N}\|\nabla g_{m}(z)\|_{L^{\infty}({\mathbb R}^{N})}.
\end{equation}
Let us use the closed ball in the function space,
\begin{equation}
\label{M}
D_{M}:=\{\nabla g(z)\in L^{\infty}({\mathbb R}^{N}, {\mathbb R}^{N^{2}}) \ | \
\|\nabla g(z)\|_{L^{\infty}({\mathbb R}^{N}, {\mathbb R}^{N^{2}})}\leq M \}, \quad M>0
\end{equation}
and impose the following conditions on the nonlinear part of the system
of equations ~\eqref{h}.

\noindent
{\bf Assumption 2.} {\it Let
$1\leq m\leq N, \ g_{m}(z): {\mathbb R}^{N}\to {\mathbb R}$, so that
$g_{m}(0)=0$. Let in addition $\nabla g(z)\in D_{M}$ and $g(z)$
is nontrivial on ${\mathbb R}^{N}$.}

Note that the solvability of a local elliptic problem in a bounded domain in
${\mathbb R}^{N}$ was studied in ~\cite{BO86}, where the nonlinear function was
allowed to have a sublinear growth.

We introduce the operator $T_g$ by $u = T_{g}v$, where $u$ is a
solution of system ~\eqref{aux}. It will be shown that, for $v\in B_{\rho}$,
problem \eqref{aux} has only one solution $u\in B_{\rho}$.  Our first main
result is as follows.

\noindent {\bf Theorem 3.} {\it Let Assumptions 1 and 2 hold and
\begin{equation}
\label{eps}
0<\varepsilon\leq \frac{\rho C}{M{\cal K}
(\|u_{0}\|_{L^{2}({\mathbb R}, {\mathbb R}^{N})}+1)}.
\end{equation}
Then the map $T_{g}: B_{\rho}\to B_{\rho}$ associated with the system of equations
~\eqref{aux} is a strict contraction.
The unique fixed point $u_{p}(x)$ of this map $T_{g}$ is the only solution of
system ~\eqref{pert} in $B_{\rho}$.}

Evidently, the resulting solution of system ~\eqref{h} given by ~\eqref{r}  will
not vanish identically on ${\mathbb R}$ since the source terms $f_{m}(x)$
are nontrivial for some $1\leq m\leq N$ and all $g_{m}(0)=0$ as we assume.

Our second main statement deals with the continuity of the cumulative solution
of the system of equations ~\eqref{h}  given by ~\eqref{r}  with respect to
the gradient of the nonlinear vector function $g$.
Let us define the following positive technical quantity
\begin{equation}
\label{sig}
\sigma:=\frac{\varepsilon M{\cal K}}{C}.
\end{equation}

\noindent
{\bf Theorem 4.} {\it Suppose that the assumptions of Theorem 3 are valid.
Let $u_{p,j}(x)$ be the unique fixed point of the map
$T_{g_{j}}: B_{\rho}\to B_{\rho}$ for $g_{j}$, $j=1$, $2$ and the cumulative
solution of system ~\eqref{h}  with $g(z)=g_{j}(z)$ be given by
\begin{equation}
\label{cum}
u_{j}(x):=u_{0}(x)+u_{p,j}(x).
\end{equation}
Then
\begin{equation}
\begin{aligned}
&\|u_{1}(x)-u_{2}(x)\|_{L^{2}({\mathbb R}, {\mathbb R}^{N})}
\\
&\phantom {u_1(x)}\leq
\frac{\varepsilon}{1-
\sigma}\frac{{\cal K}}{C}
(\|u_{0}\|_{L^{2}({\mathbb R}, {\mathbb R}^{N})}+1)\|\nabla g_{1}(z)-\nabla g_{2}(z)\|_
{L^{\infty}({\mathbb R}^{N}, {\mathbb R}^{N^{2}})}.
\end{aligned}
\label{cont}
\end{equation}
}

We turn our attention first to the proof of Theorem 3.

\bigskip

\setcounter{section}{2}
\setcounter{equation}{0}

\centerline{\bf 2. The existence of the perturbed solution}

{\it Proof of Theorem 3.} Let us choose an arbitrary $v(x)\in B_{\rho}$ and
designate the terms involved in the integral expressions in the right side of
the system of equations ~\eqref{aux}  as
\begin{equation}
\label{G}
G_{m}(x):=g_{m}(u_{0}(x)+v(x)), \quad 1\leq m\leq N.
\end{equation}
Throughout the article we will use the standard Fourier transform
\begin{equation}
\label{f}
\widehat{\phi}(p):=\frac{1}{\sqrt{2\pi}}\int_{-\infty}^{\infty}\phi(x)
e^{-ipx}dx.
\end{equation}
Clearly, the upper bound
\begin{equation}
\label{fub}
\|\widehat{\phi}(p)\|_{L^{\infty}({\mathbb R})}\leq \frac{1}{\sqrt{2\pi}}
\|\phi(x)\|_{L^{1}({\mathbb R})}
\end{equation}
holds.
We apply (\ref{f}) to both sides of system (\ref{aux}). This yields
\begin{equation}
\label{uhpKG}
\widehat{u_{m}}(p)=\varepsilon_{m} \sqrt{2\pi}
\frac{\widehat{\cal K}_{m}(p)\widehat{G_{m}}(p)}{\ln \Big(\frac{|p|}{e^{a_{m}}}
\Big)-ib_{m}p}, \quad 1\leq m\leq N.
\end{equation}
By virtue of estimates (\ref{lbab}) and (\ref{fub}), we easily obtain
\begin{equation}
\label{uhub}
|\widehat{u_{m}}(p)|\leq \varepsilon_{m} \frac{\|{\cal K}_{m}\|_{L^{1}({\mathbb R})}
|\widehat{G_{m}}(p)|}{C}, \quad 1\leq m\leq N, 
\end{equation}
such that
\begin{equation}
\label{ul2ube}
\|u_{m}\|_{L^{2}({\mathbb R})}\leq \varepsilon_{m}
\frac{\|{\cal K}_{m}\|_{L^{1}({\mathbb R})}}{C}\|G_{m}\|_{L^{2}({\mathbb R})},
\quad 1\leq m\leq N.
\end{equation}
Obviously,
\begin{equation}
\label{Gx}
G_{m}(x)=\int_{0}^{1}\nabla g_{m}(t(u_{0}(x)+v(x))).(u_{0}(x)+v(x))dt,
\quad 1\leq m\leq N,
\end{equation}
where the dot stands for the scalar product of two vectors in ${\mathbb R}^{N}$.
Hence
\begin{equation}
\label{Gaub}
|G_{m}(x)|\leq  M|u_{0}(x)+v(x)|_{{\mathbb R}^{N}}, \quad 1\leq m\leq N,
\end{equation}
where $|.|_{{\mathbb R}^{N}}$ denotes the length of a vector in ${\mathbb R}^{N}$.
Therefore, for $v(x)\in B_{\rho}$, we derive
\begin{equation}
\label{Gl2b}
\|G_{m}(x)\|_{L^{2}({\mathbb R})}\leq M
(\|u_{0}\|_{L^{2}({\mathbb R}, {\mathbb R}^{N})}+1), \quad 1\leq m\leq N.
\end{equation}
By virtue of bounds (\ref{ul2ube}) and (\ref{Gl2b}), we arrive at
\begin{equation}
\label{ul2rho}
\|u\|_{L^{2}({\mathbb R}, {\mathbb R}^{N})}\leq \frac{\varepsilon M{\cal K}}
{C}(\|u_{0}\|_{L^{2}({\mathbb R}, {\mathbb R}^{N})}+1)\leq \rho
\end{equation}
for all the values of the parameter $\varepsilon$, which satisfy ~\eqref{eps}.
This implies that $u(x)\in B_{\rho}$ as well.

Let us suppose that for some $v(x)\in B_{\rho}$ there exist two solutions
$u_{1, 2}(x)\in B_{\rho}$ of the system of equations ~\eqref{aux}. Then the
difference vector function
$w(x):=u_{1}(x)-u_{2}(x)\in L^{2}({\mathbb R},{\mathbb R}^{N} )$ solves the
homogeneous system ~\eqref{lnw}. Thus $w(x)$ vanishes a.e. on the real line
as discussed with similar arguments for Lemma 5 further down.

This means that the system of equations ~\eqref{aux}  defines a map
$T_{g}: B_{\rho}\to B_{\rho}$ for all the values of $\varepsilon$ satisfying
~\eqref{eps}.
Our goal is to establish that such map is a strict contraction.

We choose
arbitrary $v_{1,2}(x)\in B_{\rho}$. By means of the argument above
$u_{1,2}:=T_{g}v_{1,2}\in B_{\rho}$ as well when $\varepsilon$ satisfies
~\eqref{eps}.
By virtue of ~\eqref{aux}, for $1\leq m\leq N$,
$$
\left[\frac{1}{2}\ln \left(-\frac{d^{2}}{dx^{2}}\right)\right] u_{1, m}-
b_{m}\frac{du_{1, m}}{dx}-a_{m}u_{1, m} = 
$$
\begin{equation}
\label{aux1}
\varepsilon_{m} \int_{-\infty}^{\infty}{\cal K}_{m}(x-y)g_{m}(u_{0}(y)+v_{1}(y))dy,
\end{equation}
$$
\left[\frac{1}{2}\ln \left(-\frac{d^{2}}{dx^{2}}\right)\right] u_{2, m}-
b_{m}\frac{du_{2, m}}{dx}-a_{m}u_{2, m} = 
$$
\begin{equation}
\label{aux2}
\varepsilon_{m} \int_{-\infty}^{\infty}{\cal K}_{m}(x-y)g_{m}(u_{0}(y)+v_{2}(y))dy.
\end{equation}
Let us define
\begin{equation}
\label{G1G2}
G_{1, m}(x):=g_{m}(u_{0}(x)+v_{1}(x)), \quad G_{2, m}(x):=g_{m}(u_{0}(x)+v_{2}(x))
\end{equation}
with $1\leq m\leq N$.
We apply the standard Fourier transform (\ref{f}) to both sides of
systems (\ref{aux1}) and (\ref{aux2}) to obtain
\begin{equation}
\label{u12hpKG0}
\widehat{u_{1, m}}(p)=\varepsilon_{m} \sqrt{2\pi}
\frac{\widehat{\cal K}_{m}(p)\widehat{G_{1, m}}(p)}
{\ln \Big(\frac{|p|}{e^{a_{m}}}\Big)-ib_{m}p}, \quad 1\leq m\leq N,
\end{equation}    
\begin{equation}
\label{u12hpKG}
\widehat{u_{2, m}}(p)=\varepsilon_{m} \sqrt{2\pi}
\frac{\widehat{\cal K}_{m}(p)\widehat{G_{2, m}}(p)}
{\ln \Big(\frac{|p|}{e^{a_{m}}}\Big)-ib_{m}p}, \quad 1\leq m\leq N.
\end{equation}
Let us recall estimates (\ref{lbab}) and (\ref{fub}). This yields
\begin{equation}
\label{u12hub}
|\widehat{u_{1, m}}(p)-\widehat{u_{2, m}}(p)|\leq \varepsilon_{m}
\|{\cal K}_{m}\|_{L^{1}({\mathbb R})}
\frac{|\widehat{G_{1, m}}(p)-\widehat{G_{2, m}}(p)|}{C}, \quad 1\leq m\leq N.
\end{equation}
Hence
\begin{equation}
\label{u12ubn}
\|u_{1, m}-u_{2, m}\|_{L^{2}({\mathbb R})}\leq \varepsilon_{m}
\frac{\|{\cal K}_{m}\|_{L^{1}({\mathbb R})}}{C}
\|G_{1, m}-G_{2, m}\|_{L^{2}({\mathbb R})}, \quad 1\leq m\leq N.
\end{equation}
Clearly, for $1\leq m\leq N$
$$
G_{1, m}(x)-G_{2, m}(x)=
$$
\begin{equation}
\label{G12i}
\int_{0}^{1}\nabla g_{m}(u_{0}(x)+tv_{1}(x)+(1-t)v_{2}(x)).
(v_{1}(x)-v_{2}(x))dt.
\end{equation}
This allows us to estimate the norm from above as
\begin{equation}
\label{G12nub}
\|G_{1, m}-G_{2, m}\|_{L^{2}({\mathbb R})}\leq M\|v_{1}-v_{2}\|_
{L^{2}({\mathbb R}, {\mathbb R}^{N})}, \quad 1\leq m\leq N.
\end{equation}
By means of bounds (\ref{u12ubn}) and (\ref{G12nub}) along with
definitions (\ref{epsm}), (\ref{K}), (\ref{nvf}) and (\ref{sig}), we arrive at
\begin{equation}
\label{contr}
\|T_{g}v_{1}-T_{g}v_{2}\|_{L^{2}({\mathbb R}, {\mathbb R}^{N})}\leq \sigma
\|v_{1}-v_{2}\|_{L^{2}({\mathbb R}, {\mathbb R}^{N})}.
\end{equation}
By virtue of inequality (\ref{eps}) along with Remark 6 further down, we
derive
\begin{equation}
\label{sig1}
\sigma<1.
\end{equation}
Therefore, the map $T_{g}: B_{\rho}\to B_{\rho}$ defined by the system of
equations ~\eqref{aux}  is a strict contraction for all the values of
$\varepsilon$, which satisfy ~\eqref{eps}. Its unique fixed point $u_{p}(x)$ is the only
solution of system ~\eqref{pert} in the ball $B_{\rho}$. 
The cumulative 
$u(x)\in L^{2}({\mathbb R}, {\mathbb R}^{N})$ given by formula ~\eqref{r} is a
solution of problem (\ref{h}). Obviously, formula
~\eqref{ul2rho} implies that $\|u_{p}(x)\|_{L^{2}({\mathbb R}, {\mathbb R}^{N})}\to 0$ as
$\varepsilon\to 0$. \hfill\lanbox

Let us proceed to demonstrating the validity of the second main statement
of the article.
 
\bigskip


\setcounter{section}{3}
\setcounter{equation}{0}

\centerline{\bf 3. The continuity of the cumulative solution}

\bigskip

\noindent
{\it Proof of Theorem 4.} Evidently, for all the values of the parameter
$\varepsilon$ satisfying inequality ~\eqref{eps}, we have
\begin{equation}
\label{up12}
u_{p,1}=T_{g_{1}}u_{p,1}, \quad u_{p,2}=T_{g_{2}}u_{p,2}.
\end{equation}
Since 
\begin{equation}
\label{up12m}
u_{p,1}-u_{p,2}=T_{g_{1}}u_{p,1}-T_{g_{1}}u_{p,2}+T_{g_{1}}u_{p,2}-
T_{g_{2}}u_{p,2},
\end{equation}
we obtain
$$
\|u_{p,1}-u_{p,2}\|_{L^{2}({\mathbb R}, {\mathbb R}^{N})}\leq
$$
\begin{equation}
\label{up12ni}
\|T_{g_{1}}u_{p,1}-T_{g_{1}}u_{p,2}\|_{L^{2}({\mathbb R}, {\mathbb R}^{N})}+
\|T_{g_{1}}u_{p,2}-T_{g_{2}}u_{p,2}\|_{L^{2}({\mathbb R}, {\mathbb R}^{N})}.
\end{equation}
By virtue of estimate (\ref{contr}),
\begin{equation}
\label{tg1up12}
\|T_{g_{1}}u_{p,1}-T_{g_{1}}u_{p,2}\|_{L^{2}({\mathbb R}, {\mathbb R}^{N})}\leq
\sigma\|u_{p,1}-u_{p,2}\|_{L^{2}({\mathbb R}, {\mathbb R}^{N})}
\end{equation}
and bound (\ref{sig1}) is valid.
Therefore,
\begin{equation}
\label{sigma}
(1-\sigma)\|u_{p,1}-u_{p,2}\|_{L^{2}({\mathbb R}, {\mathbb R}^{N})}\leq
\|T_{g_{1}}u_{p,2}-T_{g_{2}}u_{p,2}\|_{L^{2}({\mathbb R}, {\mathbb R}^{N})}.
\end{equation}
Let $\eta(x):=T_{g_{1}}u_{p,2}$. Obviously, for $1\leq m\leq N$,
$$
\Big[\frac{1}{2}\ln \Big(-\frac{d^{2}}{dx^{2}}\Big)\Big]\eta_{m}-
b_{m}\frac{d\eta_{m}}{dx}-a_{m}\eta_{m}=
$$
\begin{equation}
\label{12}
\varepsilon_{m}
\int_{-\infty}^{\infty}{\cal K}_{m}(x-y)g_{1, m}(u_{0}(y)+u_{p,2}(y))dy,
\end{equation}
$$
\Big[\frac{1}{2}\ln \Big(-\frac{d^{2}}{dx^{2}}\Big)\Big]u_{p,2,m}-
b_{m}\frac{du_{p,2,m}}{dx}-a_{m}u_{p,2,m}=
$$
\begin{equation}
\label{22}
\varepsilon_{m}
\int_{-\infty}^{\infty}{\cal K}_{m}(x-y)g_{2,m}(u_{0}(y)+u_{p,2}(y))dy.
\end{equation}
Define
\begin{equation}
\label{G12G22}
G_{1,2,m}(x):=g_{1,m}(u_{0}(x)+u_{p,2}(x)), \quad G_{2,2,m}(x):=
g_{2,m}(u_{0}(x)+u_{p,2}(x))
\end{equation}
with $1\leq m\leq N$.
Let us apply the standard Fourier transform (\ref{f}) to both sides of
systems (\ref{12}) and (\ref{22}). This yields
\begin{equation}
\label{xiup2h}
\widehat{\eta_{m}}(p)=\varepsilon_{m} \sqrt{2 \pi}\frac{\widehat{\cal K}_{m}(p)
\widehat{G_{1,2,m}}(p)}{\ln \Big(\frac{|p|}{e^{a_{m}}}\Big)-ib_{m}p}, \quad
1\leq m\leq N,
\end{equation}
\begin{equation}
\label{xiup2h1}
\widehat{u_{p,2,m}}(p)=\varepsilon_{m} \sqrt{2 \pi}\frac{\widehat{\cal K}_{m}(p)
\widehat{G_{2,2,m}}(p)}{\ln \Big(\frac{|p|}{e^{a_{m}}}\Big)-ib_{m}p}, \quad
1\leq m\leq N.
\end{equation}
Using ~\eqref{lbab} and ~\eqref{fub}, we derive
\begin{equation}
\label{xiup2hub}
|\widehat{\eta_{m}}(p)-\widehat{u_{p,2,m}}(p)|\leq \varepsilon_{m}
\frac{\|{\cal K}_{m}\|_{L^{1}({\mathbb R})}
|\widehat{G_{1,2,m}}(p)-\widehat{G_{2,2,m}}(p)|}{C},
\end{equation}
such that
\begin{equation}
\label{xiup2hubn}
\|\eta_{m}-u_{p,2,m}\|_{L^{2}({\mathbb R})}\leq \varepsilon_{m}
\frac{\|{\cal K}_{m}\|_{L^{1}({\mathbb R})}\|G_{1,2,m}-G_{2,2,m}\|_{L^{2}({\mathbb R})}}
{C}, 
\end{equation}
where $1\leq m\leq N$. Note that
$$
G_{1,2,m}(x)-G_{2,2,m}(x)=
$$
\begin{equation}
\label{g12mint}  
\int_{0}^{1}\nabla[g_{1,m}-g_{2,m}](t(u_{0}(x)+u_{p,2}(x))).
(u_{0}(x)+u_{p,2}(x))dt, \quad 1\leq m\leq N.
\end{equation}
This enables us to obtain the estimate from above on the norm as
\begin{equation}
\label{g1222l2n}
\|G_{1,2,m}-G_{2,2,m}\|_{L^{2}({\mathbb R})}\leq
\|\nabla g_{1,m}(z)-\nabla g_{2,m}(z)\|_{L^{\infty}({\mathbb R}^{N})}
(\|u_{0}\|_{L^{2}({\mathbb R}, {\mathbb R}^{N})}+1).
\end{equation}
Let us use inequalities (\ref{sigma}), (\ref{xiup2hubn}) and
(\ref{g1222l2n}) along with definitions (\ref{epsm}) and (\ref{ng}).
Thus,
$\|u_{p,1}-u_{p,2}\|_{L^{2}({\mathbb R}, {\mathbb R}^{N})}$ can be bounded from above by
\begin{equation}
\label{up1p2h1}
\frac{\varepsilon}{1-\sigma}\frac{{\cal K}}{C}
(\|u_{0}\|_{L^{2}({\mathbb R}, {\mathbb R}^{N})}+1)
\|\nabla g_{1}(z)-\nabla g_{2}(z)\|_{L^{\infty}({\mathbb R}^{N}, {\mathbb R}^{N^{2}})}.
\end{equation}
Therefore, ~\eqref{cont} is an easy consequence of formulas ~\eqref{cum} and
~\eqref{up1p2h1}.       \hfill\lanbox

\bigskip


\setcounter{section}{4}
\setcounter{equation}{0}

\centerline{\bf 4. Auxiliary results}

Let us establish the solvability of the linear system of equations involving the
logarithmic Laplacian, the transport terms, and square-integrable right
sides
\begin{equation}
\label{lp}
\Big[\frac{1}{2}\ln \Big(-\frac{d^{2}}{dx^{2}}\Big)\Big]u_{m}-
b_{m}\frac{du_{m}}{dx}-a_{m}u_{m}=f_{m}(x), \quad
x\in {\mathbb R},
\end{equation}
where $1\leq m\leq N$ and $a_{m}, b_{m}\in {\mathbb R}, \ b_{m}\neq 0$ are
constants. The technical lemma below is the adaptation of the one used
in ~\cite{CV24} for the studies of the single integro-differential equation
involving the logarithmic Laplacian and the drift term, analogous to system
(\ref{h}). We provide it for the convenience of the readers.

\noindent
{\bf Lemma 5.} {\it  Let  $1\leq m\leq N$, the functions
$f_{m}(x):{\mathbb R}\to {\mathbb R}, \ f_{m}(x)\in L^{2} ({\mathbb R})$ and
they are nontrivial for certain $1\leq m\leq N$.  
Then the
system of equations (\ref{lp}) admits a unique solution
$u_{0}(x)\in L^{2} ({\mathbb R}, {\mathbb R}^{N})$.}

\bigskip

\noindent
{\it Proof.} First we demonstrate the uniqueness of solutions for
system (\ref{lp}). Let us suppose that it possesses two solutions
$u_{1}(x)$, $u_{2}(x)\in L^{2}({\mathbb R}, {\mathbb R}^{N})$.
Evidently, the vector function
$w(x):=u_{1}(x)-u_{2}(x)\in L^{2}({\mathbb R}, {\mathbb R}^{N})$ 
satisfies the homogeneous system of equations
\begin{equation}
\label{lnw}
\Big[\frac{1}{2}\ln \Big(-\frac{d^{2}}{dx^{2}}\Big)\Big]w_{m}-b_{m}\frac{dw_{m}}
{dx}-a_{m}w_{m}=0, \quad 1\leq m\leq N.
\end{equation}
Because each operator $L_{a_{m}, b_{m}}$ on $L^{2}({\mathbb R})$ introduced in
~\eqref{lab} has only the essential spectrum and no
nontrivial zero modes (see~\eqref{ess} and~\eqref{lbab}),
the difference of solutions $w(x)$ is trivial on ${\mathbb R}$.

Let us apply the standard Fourier transform ~\eqref{f} to both sides of
system ~\eqref{lp}. We arrive at
\begin{equation}
\label{upfp0}
\widehat{u_{m}}(p)=\frac{\widehat{f_{m}}(p)}
{\ln \Big(\frac{|p|}{e^{a_{m}}}\Big)-ib_{m}p}, \quad p\in {\mathbb R}, \quad
1\leq m\leq N.        
\end{equation}
By means of (\ref{lbab}), we have
\begin{equation}
\label{umhpub} 
|\widehat{u_{m}}(p)|\leq \frac{|\widehat{f_{m}}(p)|}{C}\in
L^{2}({\mathbb R}), \quad 1\leq m\leq N
\end{equation}
due to our assumption. Therefore, for the solution of our system of equations
~\eqref{lp}, we have $u(x)\in L^{2}({\mathbb R}, {\mathbb R}^{N})$. Its Fourier
image is given by formula ~\eqref{upfp0}.
\hfill\lanbox

\bigskip

\noindent
{\bf Remark 6.} {\it Note that under the assumptions of Lemma 5 above, the
unique solution $u_{0}(x)$ of system (\ref{lp}) belonging to
$L^{2}({\mathbb R}, {\mathbb R}^{N})$ is nontrivial because $f_{m}(x)$
do not vanish identically on the real line for some $1\leq m\leq N$.}

\bigskip


\section*{Acknowledgements}

Y.C. was supported by  NSERC of Canada
(Nos.\ RGPIN-2019-05892, RGPIN-2024-05593).

\bigskip


\section*{Declarations}

{\bf Ethical approval} The results of the work are applicable for both human
and animal studies.

\bigskip

\noindent
{\bf Publishing policy}  The authors have read and understood the publishing
policy, and submit this manuscript in accordance with this policy.

\bigskip

\noindent
{\bf Competing interests policy} The authors declare that they have
no competing interests as defined by Springer, or other interests that
might be perceived to influence the results and/or discussion reported in this
paper.

\bigskip

\noindent
{\bf Dual publication} The results in this manuscript have not been published
elsewhere, nor are they under consideration by another publisher.

\bigskip

\noindent
{\bf Authorship}  The corresponding author
has read the Springer journal policies on author responsibilities and
submits this manuscript in accordance with those policies.

\bigskip

\noindent
{\bf Author contributions statement} Yuming Chen and Vitali Vougalter 
wrote the

\noindent
manuscript. All authors reviewed the manuscript.   
      
\bigskip

\noindent
{\bf Funding} The work was partially supported by the NSERC of Canada

\noindent   
(Nos. RGPIN-2019-05892, RGPIN-2024-05593).    

\bigskip

\noindent
{\bf Third party material} All of the
material is owned by the authors and no permissions are required.

\bigskip

\noindent
{\bf Data availability statement} All data generated or analyzed during this
study are included in this published article.

\end{document}